\documentclass[aap,MSNbibl,seceqn,dvips]{arximspdf}

\doi{10.1214/12-AAP885} %kopijuoti is PTS
\volume{23}
\issue{5}
\pubyear{2013}
\firstpage{1755}
\lastpage{1777}

\makeatletter

\newcommand{\lbr}{[\![}
\newcommand{\rbr}{]\!]}

\newcommand{\sideset}[2]{\hspace*{10pt}{#1}^{\hspace*{-25pt}#2}\hspace*{10pt}}

\newtheorem{theorem}{Theorem}[section]
\newtheorem{proposition}[theorem]{Proposition}
\newtheorem{lemma}[theorem]{Lemma}
\newtheorem{corollary}[theorem]{Corollary}

\newproclaim{definition}[theorem]{Definition}
\newproclaim{remark}[theorem]{Remark}
\newproclaim{example}[theorem]{Example}
\newproclaim{assumption}[theorem]{Assumption}

\newcommand{\eps}{\varepsilon}
\newcommand{\R}{\mathbb{R}}
\newcommand{\F}{\mathbb{F}}
\renewcommand{\S}{\mathbb{S}}
\renewcommand{\L}{\mathbb{L}}
\newcommand{\cF}{\mathcal{F}}
\newcommand{\cP}{\mathcal{P}}
\newcommand{\cE}{\mathcal{E}}
\newcommand{\one}{\mathbf{1}}
\newcommand{\bD}{\mathbf{D}}

\newcommand{\Int}{\operatorname{Int}}
\newcommand{\UC}{\operatorname{UC}}
\newcommand{\esssup}{\mathop{\operatorname{ess}\operatorname{sup}}}

\newcommand{\homega}{\hat{\omega}}
\newcommand{\tomega}{\tilde{\omega}}
\newcommand{\bomega}{\bar{\omega}}

\makeatother

\begin{document}
\begin{frontmatter}

\title{Random $G$-expectations}
\runtitle{Random $G$-expectations}

\begin{aug}
\author[A]{\fnms{Marcel} \snm{Nutz}\corref{}\thanksref{t1}\ead[label=e1]{mnutz@math.columbia.edu}}
\runauthor{M. Nutz}
\affiliation{Columbia University}
\address[A]{Department of Mathematics\\
Columbia University\\
New York, New York 10027\\
USA\\
\printead{e1}} %adresu isvedimo komanda gale!
\end{aug}

\thankstext{t1}{Supported by Swiss National Science Foundation Grant
PDFM2-120424/1.}

% HISTORY:
\received{\smonth{8} \syear{2011}}
\revised{\smonth{3} \syear{2012}}

% ABSTRACT
%
\begin{abstract}
We construct a time-consistent sublinear expectation in the setting of
volatility uncertainty. This mapping extends Peng's $G$-expectation by
allowing the range of the volatility uncertainty to be stochastic. Our
construction is purely probabilistic and based on an optimal control
formulation with path-dependent control sets.
\end{abstract}

% KEYWORDS
%

\end{frontmatter}

%s1 #&#
\section{Introduction}

The so-called $G$-expectation as introduced by Peng~\cite{Peng07,Peng08} is a dynamic nonlinear expectation
which advances the notions of $g$-expectations (Peng~\cite{Peng97})
and backward SDEs (Pardoux and Peng~\cite{PardouxPeng90}). Moreover,
it yields a stochastic representation for a specific PDE and a risk
measure for volatility uncertainty in financial mathematics
(Avellaneda, Levy and Par\'as~\cite{AvellanedaLevyParas95},
Lyons~\cite{Lyons95}).
The concept of volatility uncertainty also plays a key role
in the existence theory for second order backward SDEs (Soner, Touzi
and Zhang~\cite{SonerTouziZhang2010bsde}) which were introduced as
representations for a large class of fully nonlinear second order
parabolic PDEs (Cheridito et~al.~\cite{CheriditoSonerTouziVictoir07}).\looseness=1

The $G$-expectation is a sublinear operator defined on a class of
random variables on the canonical space $\Omega$. Intuitively, it
corresponds to the ``worst-case'' expectation in a model
where the volatility of the canonical process $B$ is seen as uncertain,
but is postulated to take values in some bounded set $D$.
The symbol $G$ then stands for the support function of $D$. If $\cP^G$
is the set of martingale laws on $\Omega$ under which the volatility
of $B$ behaves accordingly,\vspace*{1pt} the $G$-expectation at time $t=0$ may be
expressed as the upper expectation
$\cE_0^G(X):=\sup_{P\in\cP^G} E^P[X]$. This description is due to
Denis, Hu and Peng~\cite{DenisHuPeng2010}. See also Denis and
Martini~\cite{DenisMartini06} for a general study of related capacities.

For positive times $t$, the $G$-expectation is extended to a
conditional expectation $\cE^G_t(X)$ with respect to the filtration
$(\cF_t)_{t\geq0}$ generated\vspace*{1pt} by $B$.
When $X=f(B_T)$ for some sufficiently regular function $f$, then $\cE^G_t(X)$ is defined via the solution of the nonlinear heat equation
$\partial_t u - G(u_{xx})=0$ with boundary condition $u|_{t=0}=f$.
The mapping $\cE^G_t$ can be extended to random variables of the form
$X=f(B_{t_1},\ldots,B_{t_n})$ by a stepwise evaluation of the PDE and
finally to a suitable completion of the space of all such random
variables. As a result, one obtains a family $(\cE_t)_{t\geq0}$ of
conditional $G$-expectations satisfying the semigroup property
$\cE_s\circ\cE_t=\cE_s$ for \mbox{$s\leq t$}, also called
time-consistency property in this context. For an exhaustive overview
of $G$-expectations and related literature we refer to Peng's recent
ICM paper~\cite{Peng10icm} and survey~\cite{Peng10}.

In this paper, we develop a formulation where the set $D$ is allowed to
be path-dependent, that is, we replace $D$ by a set-valued process $\bD
=\{\bD_t(\omega)\}$. Intuitively, this means that the function
$G(\cdot)$ is replaced by a random function $G(t,\omega,\cdot)$ and
that the a priori bounds on the volatility can be adjusted to the
observed evolution of the system, which is highly desirable for applications.
Our main result is the existence of a time-consistent family $(\cE_t)_{t\geq0}$ of sublinear operators corresponding to this formulation.
When $\bD$ depends on $\omega$ in a Markovian way, $\cE_t$ can
be seen as a stochastic representation for a class of state-dependent
nonlinear heat equations $\partial_t u -G(x,u_{xx})=0$ which are not
covered by~\cite{SonerTouziZhang2010bsde}.

At time $t=0$, we again have a set $\cP$ of probability measures and
define $\cE_0(X):=\sup_{P\in\cP} E^P[X]$. For $t>0$, we want to have
%
%e1.1 #&#
\begin{equation}
\label{eqintroDef} \cE_t(X) \mbox{\,``$=$''\,}
\sup_{P\in\cP} E^P[X|\cF_t] \qquad\mbox{in some sense.}
\end{equation}
The main difficulty here is that the set $\cP$ is not dominated by a
finite measure.
Moreover, as the resulting problem is non-Markovian in an essential
way, the PDE approach outlined above seems unfeasible. We shall adopt
the framework of regular conditional probability distributions and
define, for each $\omega\in\Omega$, a quantity $\cE_t(X)(\omega)$
by conditioning $X$ and $\bD$ (and hence, $\cP$) on the path
$\omega$ up to time $t$,
%
%e1.2 #&#
\begin{equation}
\label{eqintroDoubleDep} \cE_t(X) (\omega):=\sup_{P\in\cP(t,\omega)}
E^P\bigl[X^{t,\omega
}\bigr], \qquad\omega\in\Omega.
\end{equation}
Then the right-hand side is well defined since it is simply a supremum
of real numbers. This approach gives direct access to the underlying
measures and allows for control theoretic methods. There is no direct
reference to the function $G$, so that $G$ is no longer required to be
finite and we can work with an unbounded domain~$\bD$. The final result
is the construction of a random variable $\cE_t(X)$ which makes
(\ref{eqintroDef}) rigorous in the form
\[
\cE_t(X) = {\esssup_{P'\in\cP(t,P)}}^{(P,\cF_t)}
E^{P'}[X|\cF_t],\qquad P\mbox{-a.s.}\qquad \mbox{for all } P\in\cP,
\]
where $\cP(t,P)=\{P'\in\cP\dvtx P'=P\mbox{ on }\cF_t\}$ and $\esssup^{(P,\cF_t)}$ denotes the essential supremum with respect to the
collection of $(P,\cF_t)$-nullsets.

The approach via (\ref{eqintroDoubleDep}) is strongly inspired by the
formulation of stochastic target problems in Soner, Touzi and Zhang
\cite{SonerTouziZhang2010dual}. There, the situation is more
nonlinear in the sense that, instead of taking conditional expectations
on the right-hand side, one solves under each $P$ a backward SDE with
terminal value $X$. On the other hand, those problems have (by
assumption) a deterministic domain with respect to the volatility,
which corresponds to a deterministic set $\bD$ in our case, and
therefore their control sets are not path-dependent.

The path-dependence of $\cP(t,\omega)$
constitutes the main difficulty in the present paper, for example,
it is not obvious under which conditions $\omega\mapsto\cE_t(X)(\omega)$ in (\ref{eqintroDoubleDep}) is even measurable.
The main problem turns out to be the following.
In our formulation, the time-consistency of $(\cE_t)_{t\geq0}$ takes
the form of a dynamic programming principle. The proof of such a result
generally relies on a pasting operation performed on controls from the
various conditional problems. However, we shall see that the resulting
control in general violates the constraint given by $\bD$, when $\bD$
is stochastic.
This feature, reminiscent of viability or state-constrained problems,
renders our problem quite different from other known problems with
path-dependence, such as the controlled SDE studied by Peng
\cite{Peng04}. Our construction is based on a new notion of regularity
which is tailored such that we can perform the necessary pastings at
least on certain well-chosen controls.

One motivation for this work is to provide a model for superhedging in
financial markets with a stochastic range of volatility uncertainty.
Given a contingent claim~$X$, this is the problem of finding the
minimal capital $x$ such that by trading in the stock market $B$, one
can achieve a financial position greater or equal to $X$ at time $T$.
From a financial point of view, it is crucial that the trading strategy
be universal, that is, it should not depend on the uncertain scenario $P$.
It is worked out in Nutz and Soner~\cite{NutzSoner10} that the
(right-continuous version of the) process $\cE(X)$ yields the dynamic
superhedging price; in particular, $\cE_0(X)$ corresponds to the
minimal capital $x$. Since the universal superhedging strategy is
constructed from the quadratic covariation process of $\cE(X)$ and
$B$, it is crucial for their arguments that our model yields $\cE(X)$
as a single, aggregated process. One then obtains an ``optional
decomposition'' of the form
\[
\cE(X)=\cE_0(X) + \int Z \,dB -K,
\]
where $K$ is an increasing process whose terminal value $K_T\geq0$
indicates the difference between
the financial position time $T$ and the claim $X$.

The remainder of this paper is organized as follows. Section
\ref{sepreliminaries} introduces the basic set-up and notation.
In Section~\ref{seformulation} we formulate the control
problem (\ref{eqintroDoubleDep}) for uniformly continuous random
variables and introduce a regularity condition on $\bD$.
Section~\ref{sedynamicProg} contains the proof of the dynamic
programming principle for this control problem. In Section
\ref{seextension} we extend $\cE$ to a suitable completion.
%

%s2 #&#
\section{Preliminaries}\label{sepreliminaries}

We fix a constant $T>0$ and let $\Omega:= \{\omega\in C([0,T];\mathbb
{R}^d)\dvtx\break\omega_0=0\}$ be
the canonical space of continuous paths equipped with the uniform norm
$\|\omega\|_T:=\sup_{0\leq s\leq T} |\omega_s|$, where\vspace*{1pt}
\mbox{$|\cdot|$} is the Euclidean norm. We denote by $B$ the canonical
process $B_t(\omega)=\omega_t$, by $P_0$ the Wiener measure and by
$\F= \{\cF_t\}_{0\leq t\leq T}$ the raw filtration generated by $B$.
Unless otherwise stated, probabilistic notions requiring a filtration
(such as adaptedness) refer to $\F$.\looseness=-1

A probability measure $P$ on $\Omega$ is called \textit{local
martingale measure} if
$B$ is a local martingale under $P$.
We recall from Bichteler~\cite{Bichteler81}, Theorem 7.14, that, via
the integration-by-parts formula, the quadratic variation process
$\langle B \rangle(\omega)$ can be defined \textit{pathwise} for all
$\omega$ outside
an exceptional set which is a $P$-nullset for every local martingale
measure $P$.
Taking componentwise limits, we can then define the \mbox{$\F
$-progressively} measurable process
\[
\hat{a}_t(\omega):=\limsup_{n\to\infty} n \bigl[\langle B
\rangle_t(\omega )-\langle B \rangle_{t-1/n}(\omega) \bigr],\qquad
0<t\leq T,
\]
taking values in the set of $d\times d$-matrices with entries in the
extended real line. We also set $\hat{a}_0=0$.

Let $\bar{\cP}_W$ be the set of all local martingale measures
$P$ such that $t\mapsto\langle B \rangle_t$ is absolutely continuous
$P$-a.s.
and $\hat{a}$ takes values in $\S^{>0}_d$ $dt\times P$-a.e., where
$\S^{>0}_d\subset\R^{d\times d}$ denotes the set of strictly
positive definite matrices. Note that $\hat{a}$ is then the quadratic
variation density of $B$ under any $P\in\bar{\cP}_W$.

As in \cite
{DenisHuPeng2010,SonerTouziZhang2010dual,SonerTouziZhang2010bsde} we
shall use the so-called strong formulation of volatility uncertainty in
this paper, that is, we consider a subclass of $\bar{\cP}_W$
consisting of the laws of stochastic integrals with respect to a fixed
Brownian motion. The latter is taken to be the canonical process $B$
under $P_0$: we define $\bar{\cP}_S\subset\bar{\cP}_W$
to be the set of laws
%
%e2.1 #&#
\begin{equation}
\label{eqstrongFormulation}\quad P^\alpha:= P_0 \circ
\bigl(X^\alpha\bigr)^{-1} \qquad\mbox{where } X^{\alpha}_t:=
\sideset{\int_0^t} {(P_0)}
\alpha_s^{1/2} \,dB_s,\qquad  t\in[0,T].
\end{equation}
Here $\alpha$ ranges over all $\F$-progressively measurable processes
with values in $\S^{>0}_d$ satisfying $\int_0^T |\alpha_t|
\,dt<\infty$ $P_0$-a.s. The stochastic integral is the It\^o integral
under~$P_0$, constructed as an $\F$-progressively measurable process
with right-continuous and $P_0$-a.s. continuous paths, and, in
particular, without passing to the augmentation of $\F;$ cf. Stroock
and Varadhan~\cite{StroockVaradhan79}, page 97.

%s2.1 #&#
\subsection{Shifted paths and regular conditional distributions}

We now introduce the notation for the conditional problems of our
dynamic programming.
Since $\Omega$ is the canonical space, we can construct for any
probability measure $P$ on $\Omega$ and any $(t,\omega)\in
[0,T]\times\Omega$ the corresponding regular conditional probability
distribution~$P^\omega_t$; cf.~\cite{StroockVaradhan79},
Theorem 1.3.4. We recall that $P^\omega_t$ is a probability kernel on
$\cF_t\times\cF_T$, that is, it is a probability measure on $(\Omega
,\cF_T)$ for fixed $\omega$ and $\omega\mapsto P^\omega_t(A)$ is
$\cF_t$-measurable for each $A\in\cF_T$. Moreover, the expectation
under $P^\omega_t$ is the conditional expectation under~$P$:
\[
E^{P^\omega_t}[X]=E^P[X|\cF_t](\omega), \qquad P\mbox{-a.s.}
\]
whenever $X$ is $\cF_T$-measurable and bounded. Finally, $P^\omega_t$
is concentrated on the set of paths
that coincide with $\omega$ up to $t$,
%
%e2.2 #&#
\begin{equation}
\label{eqmeasureConcentrated} P^\omega_t \bigl\{
\omega'\in\Omega\dvtx\omega' = \omega\mbox{ on } [0,t]
\bigr\} = 1.
\end{equation}

Next, we fix $0\leq s\leq t\leq T$ and define the following shifted objects.
We denote by $\Omega^t:= \{\omega\in C([t,T];\R^d)\dvtx\omega_t=0\}$
the shifted canonical space, by $B^t$ the canonical process on
$\Omega^t$, by $P^t_0$ the Wiener measure on $\Omega^t$ and by $\F^t=\{\cF^t_u\}_{t\leq u\leq T}$ the (raw)
filtration generated by $B^t$.
For $\omega\in\Omega^s$, the shifted path $\omega^t\in\Omega^t$
is defined by
$\omega^t_u:= \omega_u-\omega_t$ for $t\leq u\leq T$ and
furthermore, if $\tomega\in\Omega^t$,
then the concatenation of $\omega$ and $\tomega$ at $t$ is the path
\[
(\omega\otimes_t \tomega)_u:= \omega_u
\one_{[s,t)}(u) + (\omega_t + \tomega_u)
\one_{[t, T]}(u),\qquad s\leq u\leq T.
\]
If $\bomega\in\Omega$, we note the associativity
$\bomega\otimes_s (\omega\otimes_t \tomega) = (\bomega\otimes_s
\omega) \otimes_t \tomega$.
Given an $\cF^s_{T}$-measurable random variable $\xi$ on $\Omega^s$
and $\omega\in\Omega^s$,
we define the shifted random variable $\xi^{t,\omega}$ on $\Omega^t$ by
\[
\xi^{t, \omega}(\tomega):=\xi(\omega\otimes_t \tomega),\qquad \tomega
\in\Omega^t.
\]
Clearly $\tomega\mapsto\xi^{t, \omega}(\tomega)$ is $\cF^t_{T}$-measurable and $\xi^{t,\omega}$ depends only on the
restriction of $\omega$ to $[s,t]$. For a random variable $\psi$ on
$\Omega$, the associativity of the concatenation yields
\[
\bigl(\psi^{s,\bomega}\bigr)^{t,\omega}=\psi^{t,\bomega\otimes_s \omega}.
\]
We note that for an $\F^s$-progressively measurable process $\{X_u,
u\in[s,T]\}$,
the shifted process $\{X^{t, \omega}_u,   u\in[t,T]\}$ is $\F^t$-progressively measurable.
If $P$ is a probability on $\Omega^s$, the measure $P^{t,\omega}$ on
$\cF^t_T$ defined by
\[
P^{t,\omega}(A):=P^\omega_t(\omega\otimes_t
A),\qquad A\in\cF^t_T\qquad \mbox{where }\omega\otimes_t
A:=\{\omega\otimes_t \tomega \dvtx\tomega\in A\}
\]
is again a probability by (\ref{eqmeasureConcentrated}). We then have
\[
E^{P^{t,\omega}}\bigl[\xi^{t,\omega}\bigr]=E^{P^\omega_t}[\xi]
=E^P\bigl[\xi|\cF^s_t\bigr](\omega) ,\qquad
P\mbox{-a.s.}
\]
In analogy to the above, we also introduce the set $\bar{\cP
}_W^t$ of martingale
measures on $\Omega^t$ under which the quadratic variation density
process $\hat{a}^t$ of $B^t$
is well defined with values in $\S^{>0}_d$ and the subset $\bar
{\cP}_S^t\subseteq\bar{\cP}_W^t$ induced
by $(P_0^t,B^t)$-stochastic integrals of $\F^t$-progressively
measurable integrands. (By convention,
$\bar{\cP}_S^T=\bar{\cP}_W^T$ consists of the unique
probability on $\Omega^T=\{0\}$.)
Finally, we denote by $\Omega^s_t:=\{\omega|_{[s,t]}\dvtx\omega\in
\Omega^s\}$ the restriction of $\Omega^s$ to $[s,t]$ and note that
$\Omega^s_t$ can be identified with $\{\omega\in\Omega^s\dvtx\omega_u=\omega_t\mbox{ for }u\in[t,T]\}$.

%s3 #&#
\section{Formulation of the control problem}\label{seformulation}

We start with a closed set-valued process $\bD\dvtx\Omega\times[0,T]\to
2^{\S_d^{+}}$ taking values in the positive semidefinite matrices,
that is, $\bD_t(\omega)$ is a closed set of matrices for each
$(t,\omega)\in[0,T]\times\Omega$. We assume that $\bD$ is
progressively measurable in the sense that for every compact $K\subset
\S_d^{+}$, the lower inverse image
$\{(t,\omega)\dvtx\bD_t(\omega)\cap K\neq\varnothing\}$ is a
progressively measurable subset of $[0,T]\times\Omega$. In
particular, the value of $\bD_t(\omega)$ depends only on the
restriction of $\omega$ to $[0,t]$.

In view of our setting with a nondominated set of probabilities, we
shall introduce topological regularity. As a first step to obtain some
stability, we consider laws under which the quadratic variation density
of $B$ takes values in a uniform interior of $\bD$.
For a set $D\subseteq\S_d^{+}$ and $\delta>0$, we define the $\delta
$-interior
$\Int^\delta D:=\{x\in D\dvtx B_\delta(x)\subseteq D\}$,
where $B_\delta(x)$ denotes the open ball of radius
$\delta$.\looseness=-1
%
%de3.1 #&#
\begin{definition}\label{defcPshiftedOpen}
Given $(t,\omega)\in[0,T]\times\Omega$, we define $\cP(t,\omega)$
to be the collection of all $P\in\bar{\cP}^t_S$ for which there
exists $\delta=\delta(t,\omega,P)>0$ such that
\[
\hat{a}^{t}_s(\tomega) \in\Int^\delta
\bD^{t,\omega}_s(\tomega ) \qquad\mbox{for } ds\times P\mbox{-a.e. } (s,
\tomega) \in [t,T]\times \Omega^t.
\]
Furthermore, if $\delta^*$ denotes the supremum of all such $\delta$,
we define the positive quantity $\deg(t,\omega,P):=(\delta^*/2)\wedge1$.
We note that $\cP(0,\omega)$ does not depend on $\omega$ and denote
this set by $\cP$.
\end{definition}

The formula $(\delta^*/2)\wedge1$ ensures that $\deg(t,\omega,P)$
is finite and among the admissible $\delta$.
The following is the main regularity condition in this paper.
%
%de3.2 #&#
\begin{definition}\label{defDunifCont}
We say that $\bD$ is \textit{uniformly continuous} if for all $\delta
>0$ and $(t,\omega)\in[0,T]\times\Omega$ there exists $\eps=\eps
(t,\omega,\delta)>0$ such that $\|\omega-\omega'\|_t\leq\eps$ implies
\[
\Int^\delta\bD^{t,\omega}_s(\tomega)\subseteq
\Int^\eps\bD^{t,\omega'}_s(\tomega) \qquad\mbox{for all }(s,
\tomega) \in [t,T]\times\Omega^t.
\]
\end{definition}

If the dimension is $d=1$ and $\bD$ is a random interval, this
property is related to the uniform continuity of the processes
delimiting the interval (see also Example~\ref{exunifContD}).
%
%as3.3 #&#
\begin{assumption}\label{asDunifcont}
We assume throughout that $\bD$ is uniformly continuous and such that
$\cP(t,\omega)\neq\varnothing$ for all $(t,\omega)\in[0,T]\times
\Omega$.
\end{assumption}

This assumption is in force for the entire paper.
We now introduce the value function which will play the role of the
sublinear (conditional) expectation. We denote by $\UC_b(\Omega)$ the
space of bounded uniformly continuous functions on $\Omega$.
%
%de3.4 #&#
\begin{definition}\label{defvalueFunction}
Given $\xi\in\UC_b(\Omega)$, we define for each $t\in[0,T]$ the
value function
\[
V_t(\omega):=V_t(\xi) (\omega):=\sup_{P\in\cP(t,\omega)}
E^P\bigl[\xi^{t,\omega}\bigr],\qquad \omega\in\Omega.\vadjust{\goodbreak}
\]
\end{definition}

Until Section~\ref{seextension}, \textit{the function $\xi$ is fixed}
and often suppressed in the notation.
The following result will guarantee enough separability for our proof
of the dynamic programming principle; it is a direct consequence of the
preceding definitions.
%
%le3.5 #&#
\begin{lemma}\label{lecontDandSeparableControls}
Let $(t,\omega)\in[0,T]\times\Omega$ and $P\in\cP(t,\omega)$.
Then there exists $\eps=\eps(t,\omega,P)>0$ such that $P\in\cP
(t,\omega')$ and $\deg(t,\omega',P)\geq\eps$ hold whenever \mbox{$\|\omega
-\omega'\|_t\leq\eps$}.
\end{lemma}
\begin{pf}
Let $\delta:=\deg(t,\omega,P)$. Then, by definition,
\[
\hat{a}^{t}_s(\tomega) \in\Int^\delta
\bD^{t,\omega}_s(\tomega ) \qquad\mbox{for } ds\times P\mbox{-a.e. } (s,
\tomega) \in [t,T]\times \Omega^t.
\]
Let $\eps=\eps(t,\omega,\delta)$ be as in Definition \ref
{defDunifCont} and $\omega'$ such that $\|\omega-\omega'\|_t\leq \eps$,
then
$\Int^\delta\bD^{t,\omega}_s(\tomega)\subseteq\Int^\eps\bD^{t,\omega'}_s(\tomega)$
by Assumption~\ref{asDunifcont} and hence,
\[
\hat{a}^{t}_s(\tomega) \in\Int^\eps
\bD^{t,\omega'}_s(\tomega ) \qquad\mbox{for } ds\times P\mbox{-a.e. } (s,
\tomega) \in [t,T]\times \Omega^t.
\]
That is, $P\in\cP(t,\omega')$ and $\deg(t,\omega',P)\geq\eps
(t,\omega,P):=(\eps/2)\wedge1$.
\end{pf}

A first consequence of the preceding lemma is the measurability of
$V_t$. We denote $\|\omega\|_t:=\sup_{0\leq s\leq t}|\omega_s|$.
%
%co3.6 #&#
\begin{corollary}\label{covalueFunctionLsc}
Let $\xi\in\UC_b(\Omega)$. The value function $\omega\mapsto
V_t(\xi)(\omega)$
is lower semicontinuous for $\|\cdot\|_t$ and, in particular, $\cF_t$-measurable.
\end{corollary}
\begin{pf}
Fix $\omega\in\Omega$ and $P\in\cP(t,\omega)$. Since $\xi$ is
uniformly continuous, there exists a modulus of continuity $\rho^{(\xi)}$,
\[
\bigl|\xi(\omega)-\xi\bigl(\omega'\bigr)\bigr|\leq\rho^{(\xi)}\bigl(
\bigl\|\omega-\omega'\bigr\|_T\bigr) \qquad\mbox{for all }\omega,
\omega'\in\Omega.
\]
It follows that for all $\tomega\in\Omega^t$,
%
%e3.1 #&#
\begin{eqnarray}
\label{eqmodulusXi} \bigl|\xi^{t,\omega}(\tomega)-\xi^{t,\omega'}(\tomega)\bigr|
&=&
\bigl|\xi(\omega\otimes_t\tomega)-\xi\bigl(\omega'
\otimes_t\tomega \bigr)\bigr|
\nonumber
\\[-2pt]
&\leq&\rho^{(\xi)}\bigl(\bigl\|\omega\otimes_t\tomega-
\omega'\otimes_t\tomega\bigr\|_T\bigr)
\\[-2pt]
&=& \rho^{(\xi)}\bigl(\bigl\|\omega-\omega'\bigr\|_t
\bigr).
\nonumber
\end{eqnarray}
Consider a sequence $(\omega^n)$ such that $\|\omega-\omega^n\|_t\to0$.
The preceding lemma shows that $P\in\cP(t,\omega^n)$ for all $n\geq
n_0=n_0(t,\omega,P)$
and thus
\begin{eqnarray*}
\liminf_{n\to\infty} V_t\bigl(\omega^n\bigr) &=&
\liminf_{n\to\infty} \sup_{P'\in\cP(t,\omega^n)} E^{P'}\bigl[
\xi^{t,\omega^n}\bigr]
\\[-2pt]
&\geq&\liminf_{n\to\infty} \Bigl[\sup_{P'\in\cP(t,\omega^n)} E^{P'}\bigl[
\xi^{t,\omega}\bigr] - \rho^{(\xi)}\bigl(\bigl\|\omega-\omega^n
\bigr\|_t\bigr) \Bigr]
\\[-2pt]
&=& \liminf_{n\to\infty} \sup_{P'\in\cP(t,\omega^n)} E^{P'}\bigl[
\xi^{t,\omega}\bigr]
\\[-2pt]
&\geq& E^P\bigl[\xi^{t,\omega}\bigr].
\end{eqnarray*}
As $P\in\cP(t,\omega)$ was arbitrary, we conclude that $\liminf_{n}
V_t(\omega^n)\geq V_t(\omega)$.\vadjust{\goodbreak}
\end{pf}

We note that the obtained regularity of $V_t$ is significantly weaker
than the uniform continuity of $\xi$; this is a consequence of the
state-dependence in our problem. Indeed, the above proof shows that if
$\cP(t,\omega)$ is independent of $\omega$, then
$V_t$ is again uniformly continuous with the same modulus of continuity
as $\xi$ (see also~\cite{SonerTouziZhang2010dual}). Similarly, in
Peng's construction of the $G$-expectation, the preservation of
Lipschitz-constants arises because the nonlinearity in the underlying
PDE has no state-dependence.%
%
%re3.7 #&#
\begin{remark}
Since $\xi$ is bounded and continuous, the value function $V_t(\xi)$
remains unchanged if $\cP(t,\omega)$ is replaced by its weak closure
(in the sense of weak convergence of probability measures). As an
application, we show that we retrieve Peng's $G$-expectation under a
nondegeneracy condition.

Given $G$, we recall from~\cite{DenisHuPeng2010}, Section 3, that
there exists a compact and convex set $D\subset\S_d^{+}$ such that
$2G$ is the support function of $D$ and such that $\cE_0^G(\psi)=\sup_{P\in\cP^G} E^P[\psi]$ for sufficiently regular $\psi$, where
\[
\cP^G:= \bigl\{P^\alpha\in\bar{\cP}_S\dvtx
\alpha_t(\omega)\in D \mbox{ for }dt\times P_0\mbox{-a.e. }
(t,\omega)\in [0,T]\times\Omega \bigr\}.
\]
We make the additional assumption that $D$ has nonempty interior $\Int
D$. In the scalar case $d=1$, this precisely rules out the trivial case
where $\cE_0^G$ is an expectation in the usual sense.

We then choose $\bD:=D$. In this deterministic situation, our
formulation boils down to
\[
\cP=\bigcup_{\delta>0} \bigl\{P^\alpha\in
\bar{\cP}_S\dvtx\alpha_t(\omega)\in\Int^\delta
D \mbox{ for }dt\times P_0\mbox{-a.e. } (t,\omega)\in[0,T]\times
\Omega \bigr\}.
\]
Clearly $\cP\subset\cP^G$, so it remains to show that $\cP$ is
dense. To this end, fix a point $\alpha^*\in\Int D$ and let $P^\alpha
\in\cP^G$, that is, $\alpha$ takes values in $D$. Then for $0<\eps
<1$, the process $\alpha^\eps:=\eps\alpha^* + (1-\eps)\alpha$
takes values in $\Int^\delta D$ for some $\delta>0$, due to the fact that
the disjoint sets $\partial D$ and $\{\eps\alpha^* + (1-\eps)x\dvtx
x\in D\}$ have positive distance by compactness.
We have $P^{\alpha^\eps}\in\cP$ and it follows by dominated
convergence for stochastic integrals that $P^{\alpha^\eps}\to
P^\alpha$ for $\eps\to0$.

While this shows that we can indeed recover the $G$-expectation, we
should mention that if one wants to treat only deterministic sets $\bD
$, one can use a much simpler construction than in this paper, and, in
particular, there is no need to use the sets $\Int^\delta\bD$ at all.
\end{remark}

Next, we give an example where our continuity assumption on $\bD$ is satisfied.
%
%ex3.8 #&#
\begin{example}\label{exunifContD}
We consider the case $d=1$. Let $a,b\dvtx[0,T]\times\Omega\to\R$
be progressively measurable processes satisfying $0\leq a<b$.
Assume that $a$ is uniformly continuous\vadjust{\goodbreak} in $\omega$, uniformly in
time, that is, that for all $\delta>0$ there exists $\eps>0$ such that
%
%e3.2 #&#
\begin{equation}
\label{equnifContEx} \bigl\|\omega-\omega'\bigr\|_T\leq\eps
\quad\mbox{implies}\quad \sup_{0\leq
s\leq
T} \bigl|a_s(\omega)-a_s\bigl(
\omega'\bigr)\bigr|\leq\delta.
\end{equation}
Assume that $b$ is uniformly continuous in the same sense. Then the
random interval
\[
\bD_t(\omega):=\bigl[a_t(\omega),b_t(\omega)
\bigr]
\]
is uniformly continuous. Indeed, given $\delta>0$, there exists $\eps'=\eps'(\delta)>0$ such that
$|a_s(\omega)-a_s(\omega')|<\delta/2$ for all $0\leq s\leq T$
whenever $\|\omega-\omega'\|_T\leq\eps'$, and the same for $b$. We
set $\eps:=\eps'\wedge\delta/2$. Then
for $\omega,\omega'$ such that $\|\omega-\omega'\|_t\leq\eps$, we
have that $\|\omega\otimes_t \tomega-\omega'\otimes_t \tomega\|_T=\|\omega-\omega'\|_t\leq\eps$ and hence,
\begin{eqnarray*}
\Int^\delta\bD^{t,\omega}_s(\tomega) &=&
\bigl[a_s(\omega\otimes_t \tomega)+\delta,
b_s(\omega \otimes_t \tomega)-\delta \bigr]
\\
&\subseteq& \bigl[a_s\bigl(\omega'
\otimes_t \tomega\bigr)+\eps, b_s\bigl(
\omega'\otimes_t \tomega\bigr)-\eps \bigr]
\\
&=& \Int^\eps\bD^{t,\omega'}_s(\tomega) \qquad\mbox{for all
}(s,\tomega) \in[t,T]\times\Omega^t.
\end{eqnarray*}

A multivariate version of the previous example runs as follows. Let
$A\dvtx\break[0,T]\times\Omega\to\S^{>0}_d$ and $r\dvtx[0,T]\times\Omega\to
[0,\infty)$ be two progressively measurable processes which are
uniformly continuous in the sense of (\ref{equnifContEx}) and define
the set-valued process
\[
\bD_t(\omega):= \bigl\{\Gamma\in\S^{>0}_d\dvtx \bigl|
\Gamma-A_t(\omega )\bigr|\leq r_t(\omega) \bigr\}.
\]
Then $\bD$ is uniformly continuous; the proof is a direct extension of
the above.
\end{example}

We close this section by a remark relating the ``random
$G$''-expectations to a class of state-dependent nonlinear heat equations.
%
%re3.9 #&#
\begin{remark}
We consider a Markovian case of Example~\ref{exunifContD}, where the
functions delimiting $\bD$ depend only on the current state. Indeed,
let $a,b\dvtx\R\to\R$ be bounded, uniformly continuous functions such
that $0\leq a\leq b$ and $b-a$ is bounded away from zero, and define
\[
\bD_t(\omega):=\bigl[a(\omega_t),b(\omega_t)
\bigr].
\]
(Of course, an additional time-dependence could also be included.)
Moreover, let $f\dvtx\R\to\R$ be a bounded, uniformly continuous
function and consider
%
%e3.3 #&#
\begin{eqnarray}
\label{eqPDE}
&\displaystyle -\partial_t u -G(x,u_{xx})=0,\qquad
u(T,\cdot)=f;&\nonumber\\[-8pt]\\[-8pt]
&\displaystyle G(x,q):=\sup_{a(x)\leq p \leq b(x)} pq/2.&\nonumber
\end{eqnarray}
We claim that the (unique, continuous) viscosity solution $u$ of
(\ref{eqPDE}) satisfies
%
%e3.4 #&#
\begin{equation}
\label{eqstochRep} u(0,x)=V_0(\xi) \qquad\mbox{for } \xi:=f(x+B_T).
\end{equation}
Indeed, by the standard Hamilton--Jacobi--Bellman theory, $u$ is the
value function of the control problem
\[
u(0,x)=\sup_\alpha E^{P_0}\bigl[f\bigl(x+ X_T^\alpha
\bigr)\bigr],\qquad X^\alpha_t=\int_0^t
\alpha^{1/2}_s \,dB_s,%
\]
where $\alpha$ varies over all positive, progressively measurable
processes satisfying
\[
\alpha_t\in\bD\bigl(X^\alpha_t\bigr), \qquad dt\times
P_0\mbox{-a.e.}
\]
For each such $\alpha$, let $P^\alpha$ be the law of $X^\alpha$,
then clearly
\[
u(0,x)=\sup_\alpha E^{P^\alpha}\bigl[f(x+ B_T)\bigr].
\]
It follows from (the proof of) Lemma~\ref{leadmissibilityByAlpha}
below that the laws $\{P^\alpha\}$ are in one-to-one correspondence
with $\cP$, if Definition~\ref{defcPshiftedOpen} is used with
$\delta=0$ (i.e., we use $\bD$ instead of its interior).
Let $G^\delta(x,q):=\sup_{a(x)+\delta\leq p \leq b(x)-\delta} pq/2$
be the nonlinearity corresponding to $\Int^\delta\bD$ and let
$u^\delta$ be the viscosity solution of the corresponding
equation (\ref{eqPDE}). Then the above yields
\[
u(0,x)\geq V_0(\xi) \geq u^\delta(0,x)
\]
for $\delta>0$ small (so that $b-a\geq2\delta$). It follows from the
comparison principle and stability of viscosity solutions that
$u^\delta(t,x)$ increases monotonically to $u(t,x)$ as $\delta
\downarrow0$; as a result, we have (\ref{eqstochRep}).
\end{remark}

%s4 #&#
\section{Dynamic programming}\label{sedynamicProg}

The main goal of this section is to prove the dynamic programming
principle for $V_t(\xi)$, which corresponds to the time-consistency
property of our sublinear expectation. For the case where $\bD$ is
deterministic and $V_t(\xi)\in\UC_b(\Omega)$, the relevant
arguments were previously given in~\cite{SonerTouziZhang2010dual}.
%

%s4.1 #&#
\subsection{Shifting and pasting of measures}

As usual, one inequality in the dynamic programming principle will be
the consequence of an invariance property of the control sets.
%
%le4.1 #&#
\begin{lemma}[(Invariance)]\label{leshiftedMeasureAdmissible}
Let $0\leq s\leq t\leq T$ and $\bomega\in\Omega$. If $P\in\cP
(s,\bomega)$, then
\[
P^{t,\omega}\in\cP(t,\bomega\otimes_s \omega) \qquad\mbox{for }P
\mbox{-a.e. } \omega\in\Omega^s.
\]
\end{lemma}
\begin{pf}
It is shown in~\cite{SonerTouziZhang2010dual}, Lemma 4.1, that
$P^{t,\omega}\in\bar{\cP}^t_S$ and that
under $P^{t,\omega}$, the quadratic variation density of $B^t$
coincides with the shift of $\hat{a}^s$:
%
%e4.1 #&#
\begin{equation}
\label{eqshiftStrongMeasure} \hat{a}^t_u(\tomega) =
\bigl(\hat{a}^s_u\bigr)^{t,\omega}(\tomega)
\qquad\mbox{for
} du \times P^{t,\omega}\mbox{-a.e. }(u, \tomega )\in[t,T]\times
\Omega^t
\end{equation}
and $P$-a.e. $\omega\in\Omega^s$. Let $\delta:=\deg(s,\bomega
,P)$, then
\[
\hat{a}^{s}_u\bigl(\omega'\bigr) \in
\Int^\delta\bD^{s,\bomega}_u\bigl(\omega'
\bigr) \qquad \mbox{for } du\times P\mbox{-a.e. } \bigl(u,\omega'\bigr) \in
[s,T]\times \Omega^s
\]
and hence,
\[
\hat{a}^{s}_u(\omega\otimes_t\tomega) \in
\Int^\delta\bD^{s,\bomega}_u(\omega\otimes_t
\tomega) \qquad\mbox{for } du\times P^{t,\omega}\mbox{-a.e. } (u,\tomega) \in[t,T]
\times\Omega^t.
\]
Now (\ref{eqshiftStrongMeasure}) shows that for $du\times P^{t,\omega
}$-a.e. $(u,\tomega)\in[t,T]\times\Omega^t$ we have
\[
\hat{a}^t_u(\tomega) =\bigl(\hat{a}^s_u
\bigr)^{t,\omega}(\tomega) = \hat {a}^s_u(\omega
\otimes_t \tomega) \in\Int^\delta\bD^{s,\bomega}_u(
\omega\otimes_t\tomega) =\Int^\delta\bD^{t,\bomega\otimes_s\omega}_u(
\tomega) %
\]
for $P$-a.e. $\omega\in\Omega^s$, that is, that $P^{t,\omega}\in
\cP(t,\bomega\otimes_s \omega)$.
\end{pf}

The dynamic programming principle is intimately related to a stability
property of the control sets under a pasting operation. More precisely,
it is necessary to collect $\eps$-optimizers from the conditional
problems over $\cP(t,\omega)$ and construct from them a control in
$\cP$ (if $s=0$). As a first step, we give a tractable criterion for
the admissibility of a control.
We recall the process $X^\alpha$ from (\ref{eqstrongFormulation})
and note that since it has continuous paths $P_0$-a.s., $X^\alpha$ can
be seen as a transformation of the canonical space under the Wiener measure.
%
%le4.2 #&#
\begin{lemma}\label{leadmissibilityByAlpha}
Let $(t,\omega)\in[0,T]\times\Omega$ and $P=P^\alpha\in\bar
{\cP}_S^t$. Then $P\in\cP(t,\omega)$ if and only if there exists
$\delta>0$ such that
\[
\alpha_s(\tomega) \in\Int^{\delta}\bD_s^{t,\omega}
\bigl(X^{\alpha
}(\tomega)\bigr) \qquad\mbox{for } ds\times P^t_0
\mbox{-a.e. } (s,\tomega)\in [t,T]\times\Omega^t.
\]
\end{lemma}
\begin{pf}
We first note that
\[
\bigl\langle B^t \bigr\rangle =\int_t^{\cdot}
\hat{a}^t_u\bigl(B^t\bigr) \,du,\qquad
P^\alpha\mbox{-a.s.} \quad\mbox{and}\quad \bigl\langle X^\alpha \bigr
\rangle = \int_t^{\cdot}\alpha_u
\bigl(B^t\bigr) \,du ,\qquad P^t_0\mbox{-a.s.}
\]
Recalling that $P^\alpha=P_0^t\circ(X^\alpha)^{-1}$, we observe that
the $P^\alpha$-distribution of
$ (B^t, \int_t^{\cdot} \hat{a}^t(B^t) \,du )$ coincides with
the $P^t_0$-distribution of
$ (X^\alpha, \int_t^{\cdot}\alpha(B^t) \,du )$.
By definition, $P^\alpha\in\cP(t,\omega)$ if and only if there
exists $\delta>0$ such that
\[
\hat{a}^{t}\bigl(B^t\bigr) \in\Int^\delta
\bD^{t,\omega}\bigl(B^t\bigr), \qquad ds\times P^\alpha
\mbox{-a.e.}\mbox{ on }[t,T]\times\Omega^t,
\]
and by the above, this is further equivalent to
\[
\alpha\bigl(B^t\bigr) \in\Int^\delta\bD^{t,\omega}
\bigl(X^\alpha\bigr) ,\qquad ds\times P_0^t\mbox{-a.e.}
\mbox{ on }[t,T]\times\Omega^t.
\]
This was the claim.
\end{pf}

To motivate\vspace*{2pt} the steps below, we first consider the admissibility of
pastings in general. We can paste given measures $P=P^\alpha\in
\bar{\cP}_S$ and $\hat{P}=P^{\hat{\alpha}}\in\bar{\cP
}^t_S$ at time $t$ to obtain a measure $\bar{P}$ on $\Omega$ and we
shall see that $\bar{P}=P^{\bar{\alpha}}$ for
\[
\bar{\alpha}_u(\omega)=\one_{[0,t)}(u)\alpha_u(
\omega) + \one_{[t,T]}(u)\hat{\alpha}_u\bigl(X^\alpha(
\omega)^t\bigr).
\]
Now assume that $P\in\cP$ and $\hat{P}\in\cP(t,\homega)$. By the
previous lemma, these constraints may be formulated as
$\alpha\in\Int^\delta\bD(X^{\alpha})$ and
$\hat{\alpha}\in\Int^\delta\bD(X^{\hat{\alpha}})^{t,\homega}$,
respectively. If $\bD$ is deterministic, we immediately see that
$\bar{\alpha }(\omega)\in\Int^\delta\bD$ for all $\omega\in\Omega$ and
therefore $\bar{P}\in\cP$. However, in the stochastic case we merely
obtain that the constraint on $\bar{\alpha}(\omega)$ is satisfied for
$\omega$ such that $X^\alpha(\omega)^t=\hat{\omega}$. Therefore, we
typically have $\bar{P}\notin\cP$.

The idea to circumvent this difficulty is that, due to the formulation
chosen in the previous section, there exists a neighborhood $B(\hat
{\omega})$ of $\hat{\omega}$ such that $\hat{P}\in\cP(t,\omega')$ for all $\omega'\in B(\hat{\omega})$. Therefore, the constraint
$\bar{\alpha}\in\Int^\delta\bD(X^{\bar{\alpha}})$ is verified
on the pre-image of $B(\hat{\omega})$ under $X^\alpha$. In the next
lemma, we exploit the separability of $\Omega$ to construct a sequence
of $\hat{P}$'s such that the corresponding neighborhoods cover the
space $\Omega$, and in Proposition~\ref{prpasting} below we shall
see how to obtain an admissible pasting from this sequence.
We denote $\|\omega\|_{[s,t]}:=\sup_{s\leq u\leq t}|\omega_u|$.
%
%le4.3 #&#
\begin{lemma}[(Separability)]\label{leseparability}
Let $0\leq s\leq t\leq T$ and $\bomega\in\Omega$. Given $\eps>0$, there
exist a sequence $(\homega^i)_{i\geq1}$ in $\Omega^s$,\vspace*{2pt} an
$\cF^s_t$-measurable partition $(E^i)_{i\geq1}$ of $\Omega^s$ and a
sequence $(P^i)_{i\geq1}$ in $\bar{\cP}_S^t$ such that:
\begin{longlist}
\item $\|\omega-\homega^i\|_{[s,t]}\leq\eps$ for all
$\omega\in E^i$,
\item $P^i\in\cP(t,\bomega\otimes_s \omega)$ for all
$\omega\in E^i$ and $\inf_{\omega\in E^i}\deg(t,\bomega\otimes_s
\omega,P^i)>0$,
\item $V_t(\bomega\otimes_s \homega^i)\leq E^{P^i}[\xi^{t,\bomega\otimes_s\homega^i}]+\eps$.
\end{longlist}
\end{lemma}
\begin{pf}
Fix $\eps>0$ and let $\homega\in\Omega^s$. By definition of
$V_t(\bomega\otimes_s\homega)$ there exists $P(\homega)\in\cP
(t,\bomega\otimes_s\homega)$ such that
\[
V_t(\homega)\leq E^{P(\homega)}\bigl[\xi^{t,\bomega\otimes_s\homega
}\bigr]+
\eps.
\]
Furthermore, by Lemma~\ref{lecontDandSeparableControls}, there exists
$\eps(\homega)=\eps(t,\bomega\otimes_s\homega,P(\homega))>0$
such that $P(\homega)\in\cP(t,\bomega\otimes_s\omega')$ and $\deg
(t,\bomega\otimes_s\omega',P(\homega))\geq\eps(\homega)$ for all
$\omega'\in\break B(\eps(\homega),\homega)\subseteq\Omega^s$. Here
$B(\eps,\homega):=\{\omega'\in\Omega^s\dvtx\|\homega-\homega'\|_{[s,t]}< \eps\}$ denotes the open $\|\cdot\|_{[s,t]}$-ball. By
replacing $\eps(\homega)$ with $\eps(\homega)\wedge\eps$ we may
assume that $\eps(\homega)\leq\eps$.

As the above holds for all $\homega\in\Omega^s$, the collection $\{
B(\eps(\homega),\homega)\dvtx\homega\in\Omega^s\}$ forms an open
cover of $\Omega^s$. Since the (quasi-)metric space $(\Omega^s,\|
\cdot\|_{[s,t]})$ is separable and therefore Lindel\"of, there exists
a countable subcover
$(B^i)_{i\geq1}$, where $B^i:=B(\eps(\homega^i),\homega^i)$. As a
$\|\cdot\|_{[s,t]}$-open set, each $B^i$ is $\cF^s_t$-measurable and
\[
E^1:=B^1,\qquad E^{i+1}:=B^{i+1}\setminus
\bigl(E^1 \cup\cdots\cup E^i\bigr),\qquad i\geq1,
\]
defines a partition of $\Omega^s$. It remains to set $P^i:=P(\homega^i)$
and note the fact that $\inf_{\omega\in E^i}\deg(t,\bomega\otimes_s\omega,P^i)
\geq\eps(\homega^i)>0$ for each $i\geq1$.
\end{pf}

For $A\in\cF^s_T$ we denote $A^{t,\omega}=\{\tomega\in\Omega^t\dvtx
\omega\otimes_t\tomega\in A\}$.
%
%pr4.4 #&#
\begin{proposition}[(Pasting)]\label{prpasting}
Let $0\leq s\leq t\leq T$, $\bomega\in\Omega$ and $P\in\cP
(s,\bomega)$.
Let $(E^i)_{0\leq i\leq N}$ be a finite $\cF^s_t$-measurable partition
of $\Omega^s$.
For $1\leq i\leq N$,\vspace*{2pt} suppose that\vadjust{\goodbreak} $P^i\in\bar{\cP}_S^t$ are such
that $P^i\in\cP(t,\bomega\otimes_s\omega)$ for all $\omega\in E^i$
and such that $\inf_{\omega\in E^i}\deg(t,\bomega\otimes_s\omega,P^i)>0$. Then
\[
\bar{P}(A):= P\bigl(A\cap E^0\bigr) + \sum
_{i=1}^N E^P \bigl[P^i
\bigl(A^{t,\omega}\bigr)\one_{E^i}(\omega) \bigr],\qquad A\in
\cF^s_T,
\]
defines an element of $\cP(s,\bomega)$. Furthermore:
\begin{longlist}
\item $\bar{P}=P$ on $\cF^s_t$,
\item $\bar{P}^{t,\omega}=P^{t,\omega}$ for $P$-a.e.
$\omega\in E^0$,
\item $\bar{P}^{t,\omega}=P^i$ for $P$-a.e. $\omega\in
E^i$ and $1\leq i\leq N$.
\end{longlist}
\end{proposition}
\begin{pf}
We first show that $\bar{P}\in\cP(s,\bomega)$. The proof that $\bar
{P}\in\bar{\cP}^s_S$ is the same as in \cite
{SonerTouziZhang2010dual}, Appendix, Proof of Claim (4.19); the
observation made there is that if $\alpha,\alpha^i$ are the $\F^s$-,
respectively, $\F^t$-progressively measurable processes such that
$P=P^\alpha$ and $P^i=P^{\alpha^i}$, then $\bar{P}=P^{\bar{\alpha}}$
for $\bar{\alpha}$ defined by
\begin{eqnarray*}
\bar{\alpha}_u(\omega)&:=&\one_{[s,t)}(u)
\alpha_u(\omega)\\
&&{} + \one_{[t,T]}(u) \Biggl[
\alpha_u(\omega)\one_{E^0}\bigl(X^\alpha(\omega)
\bigr) +\sum_{i=1}^N
\alpha^i_u\bigl(\omega^t\bigr)
\one_{E^i}\bigl(X^\alpha(\omega)\bigr) \Biggr]
\end{eqnarray*}
for $(u,\omega)\in[s,T]\times\Omega^s$.
To show that $\bar{P}\in\cP(s,\bomega)$, it remains to check that
\[
\hat{a}^s_u(\omega) \in\Int^\delta
\bD^{s,\bomega}_u(\omega ) \qquad\mbox{for }du\times\bar{P}\mbox{-a.e. }
(u,\omega)\in [s,T]\times \Omega^s
\]
for some $\delta>0$. Indeed, this is clear for $s\leq u\leq t$ since
both sides are adapted and $\bar{P}=P$ on $\cF^s_t$ by (i), which is
proved below.
In view of Lemma~\ref{leadmissibilityByAlpha}, it remains to show that
%
%e4.2 #&#
\begin{equation}
\label{eqproofPastingDelta}\quad
\bar{\alpha}_u(\omega)\in
\Int^\delta\bD^{s,\bomega}_u\bigl(X^{\bar
{\alpha}}(\omega)
\bigr) \qquad\mbox{for } du\times P^s_0\mbox{-a.e. } (u,\omega )
\in[t,T]\times\Omega^s.
\end{equation}
Let $A^i:=\{X^\alpha\in E^i\}\in\cF^s_t$ for $0\leq i\leq N$. Note
that $A^i$ is defined up to a $P^s_0$-nullset since $X^\alpha$ is
defined as an It\^o integral under $P^s_0$.
Let $\omega\in A^0$, then $X^\alpha(\omega)\in E^0$ and thus $\bar
{\alpha}_u(\omega)=\alpha_u(\omega)$ for $t\leq u\leq T$.
With $\delta^0:=\deg(s,\bomega,P)$, Lemma
\ref{leadmissibilityByAlpha} shows that
\begin{eqnarray}
\bar{\alpha}_u(\omega)=\alpha_u(\omega)\in
\Int^{\delta^0}\bD^{s,\bomega}_u\bigl(X^\alpha(\omega)
\bigr) &=&\Int^{\delta^0}\bD^{s,\bomega}_u\bigl(X^{\bar{\alpha}}(
\omega)\bigr)
\nonumber\\
&&\eqntext{\mbox{for } du\times P^s_0\mbox{-a.e. } (u,\omega)
\in[t,T]\times A^0.}
\end{eqnarray}
Next, consider $1\leq i\leq N$ and $\omega^i\in E^i$. By assumption,
$P^i\in\cP(t,\bomega\otimes_s\omega^i)$ and
\[
\deg\bigl(t,\bomega\otimes_s\omega^i,P^i
\bigr)\geq\delta^i:= \inf_{\omega
\in E^i}\deg\bigl(t,\bomega
\otimes_s\omega,P^i\bigr)>0.
\]
We set $\delta:=\min\{\delta^0,\ldots,\delta^N\}>0$, then
Lemma~\ref{leadmissibilityByAlpha} yields
\begin{eqnarray}
\alpha^i_u(\tomega) \in\Int^{\delta}
\bD_u^{t,\bomega\otimes
_s\omega^i}\bigl(X^{\alpha^i}(\tomega)\bigr) &=&
\Int^{\delta}  \bD_u^{s,\bomega}\bigl(\omega^i
\otimes_t X^{\alpha
^i}(\tomega)\bigr)
\nonumber\\
&&\eqntext{\mbox{for } du\times P^t_0\mbox{-a.e. } (u,\tomega)
\in[t,T]\times \Omega^t.}
\end{eqnarray}
Now let $\omega\in A^i$ for some $1\leq i\leq N$. Applying the
previous observation with $\omega^i:=X^\alpha(\omega)\in E^i$, we
deduce that
\begin{eqnarray}
\bar{\alpha}_u(\omega)&=&\alpha^i_u\bigl(
\omega^t\bigr)\in\Int^{\delta}\bD^{s,\bomega}_u
\bigl(X^\alpha(\omega)\otimes_t  X^{\alpha^i}\bigl(
\omega^t\bigr)\bigr) = \Int^{\delta} \bD^{s,\bomega}_u
\bigl(X^{\bar{\alpha}}(\omega )\bigr)
\nonumber\\
&&\eqntext{\mbox{for }du\times P^s_0\mbox{-a.e. } (u,\omega)
\in[t,T]\times A^i.}
\end{eqnarray}
More precisely, here we have used the following two facts. First,
to pass from $du\times P^t_0$-nullsets to $du\times P^s_0$-nullsets, we
have used that
if $G\subset\Omega^t$ is a $P^t_0$-nullset, then $P^s_0\{\omega\in
\Omega^s\dvtx\omega^t\in G\}=P_0^t(G)=0$ since the canonical process
$B^s$ has $P^s_0$-independent increments. Second,
we have used that $\psi(\omega):=X^\alpha(\omega)\otimes_t
X^{\alpha^i}(\omega^t)=X^{\bar{\alpha}}(\omega)$ for $\omega\in
A^i$. Indeed, for $s\leq u <t$ we have
$\psi_u(\omega)=X_u^\alpha(\omega)=X_u^{\bar{\alpha}}(\omega)$,
while for $t\leq u\leq T$, $\psi_u(\omega)$ equals
\[
\sideset{\int_s^t}{(P^s_0)}
\alpha^{1/2} \,dB(\omega) + \sideset{\int_t^u}{(P^t_0)} \bigl(\alpha^i
\bigr)^{1/2} \,dB^t \bigl(\omega^t\bigr) =\sideset{\int_s^u}{(P^s_0)}
(\bar{\alpha})^{1/2} \,dB(\omega) =X_u^{\bar{\alpha}}(
\omega).
\]
As $P^s_0 [\bigcup_{i=0}^N A^i ]=1$, we have proved
(\ref{eqproofPastingDelta}), therefore $\bar{P}\in\cP(s,\bomega)$.

It remains to show (i)--(iii). These assertions are fairly standard; we
include the proofs for completeness.

(i) Let $A\in\cF^s_t$; we show that $\bar{P}(A)=P(A)$. Indeed, for
$\omega\in\Omega$, the question whether $\omega\in A$ depends only
on the restriction of $\omega$ to $[s,t]$. Therefore,
\[
P^i\bigl(A^{t,\omega}\bigr)=P^i\{\tomega\dvtx\omega
\otimes_t \tomega\in A\}=\one_A(\omega),\qquad 1\leq i\leq N,
\]
and thus $\bar{P}(A)= \sum_{i=0}^N E^P[\one_{A\cap
E^i}]=P(A)$.\vspace*{1pt}

(ii), (iii) Let $F\in\cF^t_T$; we show that
\[
\bar{P}^{t,\omega}(F)=P^{t,\omega}(F)\one_{E^0}(\omega)+\sum
_{i=1}^N P^i(F)
\one_{E^i}(\omega) ,\qquad P\mbox{-a.s.}
\]
Using the definition of conditional expectation and (i), this is
equivalent to the following equality for all $\Lambda\in\cF^s_t$:
\[
\bar{P}\bigl\{\omega\in\Lambda\dvtx\omega^t\in F\bigr\}=P\bigl\{\omega
\in\Lambda \cap E^0\dvtx\omega^t\in F\bigr\}+ \sum
_{i=1}^N P^i(F)P\bigl(\Lambda\cap
E^i\bigr).
\]
For $A:=\{\omega\in\Lambda\dvtx\omega^t\in F\}$ we have $A^{t,\omega
}=\{\tomega\in F\dvtx\omega\otimes_t\tomega\in\Lambda\}$ and since
$\Lambda\in\cF^s_t$, $A^{t,\omega}$ equals $F$ if $\omega\in
\Lambda$ and is empty otherwise. Thus the definition of $\bar{P}$ yields
$\bar{P}(A)=P(A\cap E^0)+ \sum_{i=1}^N E^P[P^i(F)\one_{\Lambda}(\omega
)\one_{E^i}(\omega)]=P(A\cap E^0)+ \sum_{i=1}^N P^i(F)P(\Lambda\cap E^i)$
as desired.
\end{pf}

We remark that the above arguments apply also to a countably infinite
partition $(E^i)_{i\geq1}$, provided that
$\inf_{i\geq1}\inf_{\omega\in E^i}\deg(t,\omega,P^i)>0$. However,
this condition is difficult to guarantee. A second observation is that
the results of this subsection are based on the regularity property of
$\omega\mapsto\cP(t,\omega)$ stated in Lemma
\ref{lecontDandSeparableControls}, but make no use of the continuity of
$\xi$ or the measurability of $V_t(\xi)$.

%s4.2 #&#
\subsection{Dynamic programming principle}\label{subsecDPP}

We can now prove the key result of this paper. We recall the value
function $V_t=V_t(\xi)$ from Definition~\ref{defvalueFunction} and
denote by $\esssup^{(P,\cF_s)}$ the essential supremum of a family of
$\cF_s$-measurable random variables with respect to the collection of
$(P,\cF_s)$-nullsets.
%
%th4.5 #&#
\begin{theorem}\label{thDPP}
Let $0\leq s\leq t\leq T$. Then
%
%e4.3 #&#
\begin{equation}
\label{eqDPP} V_s(\omega) = \sup_{P\in\cP(s,\omega)} E^{P}
\bigl[ (V_t)^{s,\omega
} \bigr] \qquad\mbox{for all } \omega\in\Omega.
\end{equation}
With $\cP(s,P):=\{P'\in\cP\dvtx P'=P\mbox{ on }\cF_s\}$, we also have
%
%e4.4 #&#
\begin{equation}
\label{eqDPPesssup} V_s = {\esssup_{P'\in\cP(s,P)}}^{(P,\cF_s)}
E^{P'}[V_t|\cF_s] ,\qquad P\mbox{-a.s.} \qquad\mbox{for all }
P\in\cP
\end{equation}
and, in particular,
%
%e4.5 #&#
\begin{equation}
\label{eqDPPesssupSpecial} V_s = {\esssup_{P'\in\cP(s,P)}}^{(P,\cF_s)}
E^{P'}[\xi|\cF_s] ,\qquad P\mbox{-a.s.}\qquad \mbox{for all } P\in\cP.
\end{equation}
\end{theorem}
\begin{pf}
(i) We first show the inequality ``$\leq$'' in (\ref{eqDPP}).
Fix $\bomega\in\Omega$ as well as $P\in\cP(s,\bomega)$.
Lemma~\ref{leshiftedMeasureAdmissible} shows that $P^{t,\omega}\in
\cP(t,\bomega\otimes_s \omega)$ for $P$-a.e. $\omega\in\Omega^s$, yielding the inequality in
\begin{eqnarray*}
E^{P^{t,\omega}} \bigl[\bigl(\xi^{s,\bomega}\bigr)^{t,\omega} \bigr]
&=&
E^{P^{t,\omega}} \bigl[\xi^{t,\bomega\otimes_s\omega} \bigr]
\\
&\leq& \sup_{P'\in\cP(t,\bomega\otimes_s \omega)} E^{P'} \bigl[\xi^{t,\bomega\otimes_s\omega} \bigr]
\\
& = & V_t(\bomega\otimes_s \omega)
\\
& = & V_t^{s,\bomega}(\omega) \qquad\mbox{for $P$-a.e. $\omega\in
\Omega^s$},
\end{eqnarray*}
where $V_t^{s,\bomega}:=(V_t)^{s,\bomega}$. Since $V_t$ is measurable
by Corollary~\ref{covalueFunctionLsc}, we can take $P(d\omega
)$-expectations on both sides to obtain that
\[
E^P \bigl[\xi^{s,\bomega} \bigr] %
= E^{P}
\bigl[E^{P^{t,\omega}} \bigl[\bigl(\xi^{s,\bomega}\bigr)^{t,\omega} \bigr]
\bigr] \leq E^P \bigl[V_t^{s,\bomega} \bigr].
\]
Thus taking supremum over $P\in\cP(s,\bomega)$ yields the claim.

(ii) We now show the inequality ``$\geq$'' in (\ref{eqDPP}). Fix
$\bomega\in\Omega$ and $P\in\cP(s,\bomega)$ and let $\delta>0$.
We start with a preparatory step.

(ii.a) We claim that there exists a $\|\cdot\|_{[s,t]}$-compact set
$E\in\cF^s_t$ satisfying $P(E)> 1-\delta$ such that the restriction
\[
V_t^{s,\bomega}(\cdot)|_{E} \qquad\mbox{is uniformly
continuous for } \|\cdot\|_{[s,t]}.
\]
In particular, there exists a modulus of continuity $\rho^{(V^{s,\bomega}_t|E)}$ such that
\[
\bigl|V^{s,\bomega}_t(\omega)-V^{s,\bomega}_t\bigl(
\omega'\bigr)\bigr|\leq\rho^{(V^{s,\bomega}_t|E)} \bigl(\bigl\|\omega-
\omega'\bigr\|_{[s,t]} \bigr) \qquad\mbox{for all }\omega,
\omega'\in E.
\]
Indeed, since $P$ is a Borel measure on the Polish space $\Omega^s_t$,
there exists a compact set $K=K(P,\delta)\subset\Omega^s_t$ such that\vadjust{\goodbreak}
$P(K)>1-\delta/2$. As $V_t^{s,\bomega}$ is $\cF^s_t$-measurable (and
thus Borel-measurable as a function on $\Omega^s_t$), there exists, by
Lusin's theorem, a~closed set $\Lambda=\Lambda(P,\delta)\subseteq
\Omega^s_t$ such that $P(\Lambda)>1-\delta/2$ and such that
$V_t^{s,\bomega}|_{\Lambda}$ is $\|\cdot\|_{[s,t]}$-continuous. Then
$E':=K\cap\Lambda\subset\Omega^s_t$ is compact and hence, the
restriction of $V_t^{s,\bomega}$ to $E'$ is even uniformly continuous.
It remains to set $E:=\{\omega\in\Omega^s\dvtx\omega|_{[s,t]}\in E'\}$.

(ii.b) Let $\eps>0$. We apply Lemma~\ref{leseparability} to $E$
(instead of $\Omega^s$) and obtain
a sequence $(\homega^i)$ in $E$, an $\cF^s_t$-measurable partition
$(E^i)$ of $E$ and a sequence
$(P^i)$ in $\bar{\cP}_S^t$ such that:
\begin{longlist}[(c)]
\item[(a)] $\|\omega-\homega^i\|_{[s,t]}\leq\eps$ for all $\omega
\in E^i$,
\item[(b)] $P^i\in\cP(t,\bomega\otimes_s \omega)$ for all $\omega
\in E^i$ and $\inf_{\omega\in E^i}\deg(t,\bomega\otimes_s \omega,P^i)>0$,
\item[(c)] $V_t(\bomega\otimes_s \homega^i)\leq E^{P^i}[\xi^{t,\bomega\otimes_s\homega^i}]+\eps$.
\end{longlist}
Let $A_N:= E^1\cup\cdots\cup E^N$ for $N\geq1$.
In view of (a)--(c), we can apply Proposition~\ref{prpasting} to the
finite partition $\{A_N^c,E^1,\ldots,E^N\}$ of $\Omega^s$ and obtain a
measure $\bar{P}=\bar{P}_N\in\cP(s,\bomega)$ such that
\[
\bar{P}=P \qquad\mbox{on }\cF^s_t \quad\mbox{and}\quad
\bar{P}^{t,\omega}= \cases{ P^{t,\omega}, &\quad for $\omega\in
A_N^c$,
\vspace*{1pt}\cr
P^i, &\quad for $\omega\in
E^i, 1\leq i\leq N$.}
\]
Since $\xi$ is uniformly continuous, we obtain, similar to (\ref
{eqmodulusXi}), that there exists a modulus of continuity
$\rho^{(\xi)}$ such that
\[
\bigl|\xi^{t,\bomega\otimes_s\omega}-\xi^{t,\bomega\otimes_s\omega
'}\bigr|\leq\rho^{(\xi)}\bigl(\bigl\|\omega-
\omega'\bigr\|_{[s,t]}\bigr).
\]
Let $\omega\in E^i\subset\Omega^s$ for some $1\leq i\leq N$. Then
using (a) and (c),
\begin{eqnarray*}
V^{s,\bomega}_t(\omega) & \leq & V^{s,\bomega}_t
\bigl(\homega^i\bigr) + \rho^{(V^{s,\bomega}_t|E)}(\eps )
\\
& \leq & E^{P^i} \bigl[\xi^{t,\bomega\otimes_s\homega^i} \bigr] + \eps +
\rho^{(V^{s,\bomega}_t|E)}(\eps)
\\
& \leq & E^{P^i} \bigl[\xi^{t,\bomega\otimes_s\omega} \bigr] + \rho^{(\xi)}(
\eps) + \eps+ \rho^{(V^{s,\bomega}_t|E)}(\eps)
\\
& = & E^{\bar{P}^{t,\omega}} \bigl[\xi^{t,\bomega\otimes_s\omega
} \bigr] + \rho^{(\xi)}(\eps)
+ \eps+ \rho^{(V^{s,\bomega
}_t|E)}(\eps)
\\
& = & E^{\bar{P}^{t,\omega}} \bigl[\bigl(\xi^{s,\bomega}\bigr)^{t,\omega} \bigr]
+ \rho^{(\xi)}(\eps) + \eps+ \rho^{(V^{s,\bomega}_t|E)}(\eps)
\\
& = & E^{\bar{P}} \bigl[\xi^{s,\bomega} |\cF^s_t
\bigr](\omega) + \rho^{(\xi)}(\eps) + \eps+ \rho^{(V^{s,\bomega}_t|E)}(\eps)
\end{eqnarray*}
for $\bar{P}$-a.e. (and thus $P$-a.e.) $\omega\in E^i$. This holds
for all $1\leq i\leq N$.
As $P=\bar{P}$ on~$\cF^s_t$, taking $P$-expectations yields
\[
E^P\bigl[V^{s,\bomega}_t \one_{A_N}\bigr]
\leq E^{\bar{P}}\bigl[\xi^{s,\bomega} \one_{A_N}\bigr] +
\rho^{(\xi)}(\eps) + \eps+ \rho^{(V^{s,\bomega
}_t|E)}(\eps).
\]
Recall that $\bar{P}=\bar{P}_N$. Using dominated convergence on the
left-hand side, and on the right-hand side that $\bar{P}_N(E\setminus
A_N)=P(E\setminus A_N)\to0$ as $N\to\infty$ and that
%
%e4.6 #&#
\begin{eqnarray}
\label{eqproofDPPlimit} E^{\bar{P}_N}\bigl[\xi^{s,\bomega}
\one_{A_N}\bigr] &=& E^{\bar{P}_N}\bigl[\xi^{s,\bomega}
\one_{E}\bigr] - E^{\bar{P}_N}\bigl[\xi^{s,\bomega}
\one_{E\setminus A_N}\bigr]
\nonumber\\[-8pt]\\[-8pt]
&\leq& E^{\bar{P}_N}\bigl[\xi^{s,\bomega} \one_{E}\bigr] + \|
\xi\|_{\infty} P_N(E\setminus A_N),\nonumber
\end{eqnarray}
we conclude that
\begin{eqnarray*}
E^P\bigl[V^{s,\bomega}_t \one_E\bigr]
&\leq& \limsup_{N\to\infty}E^{\bar{P}_N}\bigl[\xi^{s,\bomega}
\one_{E}\bigr] + \rho^{(\xi)}(\eps) + \eps+
\rho^{(V^{s,\bomega}_t|E)}(\eps)
\\
&\leq& \sup_{P'\in\cP(s,\bomega,t,P)} E^{P'}\bigl[\xi^{s,\bomega}
\one_E\bigr]+ \rho^{(\xi)}(\eps) + \eps+ \rho^{(V^{s,\bomega}_t|E)}(
\eps),
\end{eqnarray*}
where $\cP(s,\bomega,t,P):=\{P'\in\cP(s,\bomega)\dvtx P'=P\mbox{ on
}\cF^s_t\}$.
As $\eps>0$ was arbitrary, this shows that
\[
E^P\bigl[V^{s,\bomega}_t \one_E\bigr]
\leq\sup_{P'\in\cP(s,\bomega,t,P)} E^{P'}\bigl[\xi^{s,\bomega}
\one_E\bigr].
\]
Finally, since $P'(E)=P(E)>1-\delta$ for all $P'\in\cP(s,\bomega
,t,P)$ and $\delta>0$ was arbitrary, we obtain by an argument similar
to (\ref{eqproofDPPlimit}) that
\[
E^P\bigl[V^{s,\bomega}_t\bigr] \leq
\sup_{P'\in\cP(s,\bomega,t,P)} E^{P'}\bigl[\xi^{s,\bomega}\bigr] \leq
\sup_{P'\in\cP(s,\bomega)} E^{P'}\bigl[\xi^{s,\bomega}\bigr] =
V_s(\bomega).
\]
The claim follows as $P\in\cP(s,\bomega)$ was arbitrary. The proof
of (\ref{eqDPP}) is complete.

(iii) The next step is to prove that
%
%e4.7 #&#
\begin{equation}
\label{eqDPPesssupIntermediate} V_t\leq {\esssup_{P'\in\cP(t,P)}}^{(P,\cF
_t)}
E^{P'}[\xi|\cF_t] ,\qquad P\mbox{-a.s.} \qquad\mbox{for all } P\in \cP.
\end{equation}
Fix $P\in\cP$. We use the previous step (ii) for the special case
$s=0$ and obtain that given $\eps>0$ there exists for each $N\geq
1$ a measure $\bar{P}_N\in\cP(t,P)$ such that
\begin{eqnarray}
V_t(\omega)&\leq& E^{\bar{P}_N}[\xi| \cF_t](\omega) +
\rho^{(\xi
)}(\eps) + \eps+ \rho^{(V_t|E)}(\eps) \nonumber\\
&&\eqntext{\mbox{for }P\mbox
{-a.s. } \omega\in E^1\cup\cdots\cup E^N.}
\end{eqnarray}
Therefore, since $E=\bigcup_{i\geq1} E^i$, we have
\[
V_t(\omega)\leq\sup_{N\geq1}E^{\bar{P}_N}[\xi|
\cF_t](\omega) + \rho^{(\xi)}(\eps) + \eps+
\rho^{(V_t|E)}(\eps) \qquad\mbox{for }P\mbox{-a.s. } \omega\in E.
\]
We recall that the set $E$ depends on $\delta$, but not on $\eps$.
Thus, letting $\eps\to0$ yields
\begin{eqnarray*}
V_t\one_E&\leq& {\esssup_{P'\in\cP(t,P)}}^{
(P,\cF_t)}
\bigl(E^{P'}[\xi| \cF_t]\one_E \bigr)\\
&=& \Bigl(
{\esssup_{P'\in\cP(t,P)}}^{(P,\cF_t)}E^{P'}[\xi| \cF_t]
\Bigr)\one_E ,\qquad P\mbox{-a.s.},
\end{eqnarray*}
where we have used that $E\in\cF_t$. In view of $P(E)> 1-\delta$,
the claim follows by taking the limit $\delta\to0$.

\mbox{}\hspace*{1pt}(iv) We now prove the inequality ``$\leq$'' in (\ref{eqDPPesssup});
we shall reduce this claim to its special case (\ref
{eqDPPesssupIntermediate}). Fix $P\in\cP$. For any $P'\in\cP(s,P)$
we have that $(P')^{t,\omega}\in\cP(t,\omega)$ for $P'$-a.s.
$\omega\in\Omega$ by Lemma~\ref{leshiftedMeasureAdmissible}.\vadjust{\goodbreak}
Thus we can infer from (\ref{eqDPP}), applied with $s:=t$ and $t:=T$, that
\[
V_t(\omega) \geq E^{(P')^{t,\omega}}\bigl[
\xi^{t,\omega}\bigr]=E^{P'}[\xi |\cF_t](\omega),
\qquad P'\mbox{-a.s.}
\]
and, in particular, that $E^{P'}[V_t|\cF_s]\geq E^{P'}[\xi|\cF_s]$
$P'$-a.s. on $\cF_s$, hence, also $P$-a.s. This shows that
\[
{\esssup_{P'\in\cP(s,P)}}^{(P,\cF_s)} E^{P'}[V_t|
\cF_s] \geq {\esssup_{P'\in\cP(s,P)}}^{
(P,\cF_s)} E^{P'}[
\xi|\cF_s] ,\qquad P\mbox{-a.s.}
\]
But (\ref{eqDPPesssupIntermediate}), applied with $s$ instead of $t$,
yields that the right-hand side $P$-a.s. dominates $V_s$. This proves
the claim.

(v) It remains to show the inequality ``$\geq$'' in (\ref
{eqDPPesssup}). Fix $P\in\cP$ and $P'\in\cP(s,P)$. Since
$(P')^{s,\omega}\in\cP(s,\omega)$ for $P'$-a.s. $\omega\in\Omega$
by Lemma~\ref{leshiftedMeasureAdmissible}, (\ref{eqDPP}) yields
\[
V_s(\omega) \geq E^{(P')^{s,\omega}}\bigl[V_t^{s,\omega}
\bigr]=E^{P'}[V_t|\cF_s](\omega)
\]
$P'$-a.s. on $\cF_s$ and hence, also $P$-a.s. The claim follows as
$P'\in\cP(s,P)$ was arbitrary.
\end{pf}

%s5 #&#
\section{Extension to the completion}\label{seextension}

So far, we have studied the value function $V_t=V_t(\xi)$ for $\xi\in
\UC_b(\Omega)$. We now set $\cE_t(\xi):=V_t$ and extend this
operator to a completion of $\UC_b(\Omega)$ by the usual procedure
(e.g., Peng~\cite{Peng05}). The main result in this section is that
the dynamic programming principle carries over to the extension.

Given $p\in[1,\infty)$ and $t\in[0,T]$, we define $L^p_\cP(\cF_t)$
to be the space of $\cF_t$-measurable random variables $X$ satisfying
\[
\|X\|_{L^p_\cP}:=\sup_{P\in\cP} \|X\|_{L^p(P)}<\infty,
\]
where $\|X\|^p_{L^p(P)}:=E^P[|X|^p]$. More precisely, we take
equivalences classes with respect to $\cP$-quasi-sure equality so that
$L^p_\cP(\cF_t)$ becomes a Banach space. [Two functions are equal
$\cP$-quasi-surely ($\cP$-q.s. for short) if they are equal
$P$-a.s. for all $P\in\cP$.] Furthermore,
\[
\L^p_\cP(\cF_t)\mbox{ is defined as the $\|
\cdot\|_{L^p_\cP
}$-closure of }\UC_b(\Omega_t)
\subseteq L^p_\cP(\cF_t).
\]
For brevity, we shall sometimes write $\L^p_\cP$ for $\L^p_\cP(\cF_T)$ and $L^p_\cP$ for $L^p_\cP(\cF_T)$.
%
%le5.1 #&#
\begin{lemma}\label{leLipschitzAndExtension}
Let $p\in[1,\infty)$. The mapping $\cE_t$ on $\UC_b(\Omega)$ is
$1$-Lipschitz for the norm $\|\cdot\|_{L^p_\cP}$,
\[
\bigl\|\cE_t(\xi)-\cE_t(\psi)\bigr\|_{L^p_\cP}\leq\|\xi-\psi
\|_{L^p_\cP
} \qquad\mbox{for all } \xi,\psi\in\UC_b(\Omega).
\]
As a consequence, $\cE_t$ uniquely extends to a Lipschitz-continuous mapping
\[
\cE_t\dvtx\L^p_\cP(\cF_T)\to
L^p_\cP(\cF_t).
\]
\end{lemma}
\begin{pf}
Note that $|\xi-\psi|^p$ is again in $\UC_b(\Omega)$. The
definition of $\cE_t$ and Jensen's inequality imply that\vadjust{\goodbreak}
$|\cE_t(\xi)-\cE_t(\psi)|^p\leq\cE_t(|\xi-\psi|)^p\leq\cE_t(|\xi-\psi|^p)$. Therefore,
\[
\bigl\|\cE_t(\xi)-\cE_t(\psi)\bigr\|_{L^p_\cP} \leq
\sup_{P\in\cP} E^P \bigl[\cE_t\bigl(|\xi-
\psi|^p\bigr) \bigr]^{1/p} = \sup_{P\in\cP}
E^P\bigl[ |\xi-\psi|^p \bigr]^{1/p},
\]
where the equality is due to (\ref{eqDPP}).
\end{pf}

Since we shall use $\L^p_\cP$ as the domain of $\cE_t$, we also give
an explicit description of this space.
We say that (an equivalence class) $X\in L^1_\cP$ is
\textit{$\cP$-quasi uniformly continuous} if $X$ has a representative $X'$ with the
property that for all $\eps>0$ there exists an open set $G\subset
\Omega$ such that $P(G)<\eps$ for all $P\in\cP$ and such that the
restriction $X'|_{\Omega\setminus G}$ is uniformly continuous. We
define $\cP$-quasi continuity in an analogous way and denote by
$C_b(\Omega)$ the space of bounded continuous functions on $\Omega$.
The following is very similar to the results in~\cite{DenisHuPeng2010}.

%pr5.2 #&#
\begin{proposition}\label{prdescriptionCompletion}
Let $p\in[1,\infty)$. The space $\L^p_\cP$ consists of all $X\in
L^p_\cP$ such that $X$ is $\cP$-quasi uniformly continuous and $\lim_{n\to\infty} \|X\one_{\{|X|\geq n\}}\|_{L^p_\cP}=0$.

If $\bD$ is uniformly bounded, then $\L^p_\cP$ coincides with the \mbox{$\|
\cdot\|_{L^p_\cP}$}-closure of $C_b(\Omega)\subset L^p_\cP$ and
``uniformly continuous'' can be replaced by ``continuous.''
\end{proposition}
\begin{pf}
For the first part, it suffices to go through the proof of Theorem~25
of~\cite{DenisHuPeng2010} and replace continuity by uniform
continuity everywhere. The only difference is that one has to use a
refined version of Tietze's extension theorem which yields uniformly
continuous extensions; cf. Mandelkern~\cite{Mandelkern90}.

If $\bD$ is uniformly bounded, $\cP$ is a set of laws of continuous
martingales with uniformly bounded quadratic variation density and
therefore $\cP$ is tight. Together with the aforementioned extension
theorem we derive that $C_b(\Omega)$ is contained in $\L^p_\cP$ and
now the result follows from~\cite{DenisHuPeng2010}, Theorem 25.
\end{pf}

Before extending the dynamic programming principle, we prove the
following auxiliary result which shows, in particular, that the
essential suprema in Theorem~\ref{thDPP} can be represented as
increasing limits.
This is a consequence of a standard pasting argument which involves
only controls with the same ``history'' and hence, there are no
problems of admissibility as in Section~\ref{sedynamicProg}.
%
%le5.3 #&#
\begin{lemma}\label{leincreasingSequence}
Let $\tau$ be an $\F$-stopping time and $X\in L^1_\cP(\cF_T)$. For
each $P\in\cP$ there exists a sequence $P_n\in\cP(\tau,P)$ such that
\[
{\esssup_{P'\in\cP(\tau,P)}}^{(P,\cF_\tau
)} E^{P'}[X|\cF_\tau]=
\lim_{n\to\infty} E^{P_n}[X|\cF_\tau] ,\qquad P\mbox{-a.s.},
\]
where the limit is increasing and $\cP(\tau,P):=\{P'\in\cP\dvtx
P'=P\mbox{ on }\cF_\tau\}$.
\end{lemma}
\begin{pf}
It suffices to show that the set $\{E^{P'}[X|\cF_\tau]\dvtx P'\in\cP
(\tau,P)\}$ is
$P$-a.s. upward filtering. Indeed,\vadjust{\goodbreak} we prove that for $\Lambda\in\cF_\tau$ and $P_1,P_2\in\cP(\tau,P)$
there exists $\bar{P}\in\cP(\tau,P)$ such that
\[
E^{\bar{P}}[X|\cF_\tau]=E^{P_1}[X|\cF_\tau]
\one_\Lambda+ E^{P_2}[X|\cF_\tau]\one_{\Lambda^c},\qquad
P\mbox{-a.s.},
\]
then the claim follows by letting $\Lambda:=\{E^{P_1}[X|\cF_\tau
]>E^{P_2}[X|\cF_\tau]\}$.
Similarly as in Proposition~\ref{prpasting}, we define
%
%e5.1 #&#
\begin{equation}
\label{eqproofPastingStop} \bar{P}(A):= E^P \bigl[P^1(A|
\cF_\tau)\one_{\Lambda} + P^2(A|\cF_\tau
)\one_{\Lambda^c} \bigr],\qquad A\in\cF_T.
\end{equation}
Let $\alpha,\alpha^1,\alpha^2$ be such that $P^{\alpha}=P$,
$P^{\alpha^1}=P_1$ and $P^{\alpha^2}=P_2$. The fact that $P=P^1=P^2$
on $\cF_\tau$ translates to $\alpha=\alpha^1=\alpha^2$ $du\times
P_0$-a.e. on $\lbr0,\tau(X^\alpha)\lbr$
and with this observation we have as in Proposition~\ref{prpasting}
that $\bar{P}=P^{\bar{\alpha}}\in\bar{\cP}_S$ for
the $\F$-progressively measurable process
\begin{eqnarray*}
\bar{\alpha}_u(\omega)&:=&
\one_{\lbr0,\tau(X^\alpha)\lbr}(u)\alpha_u(\omega) \\
&&{}+ \one_{\lbr\tau
(X^\alpha),T\rbr}(u)
\bigl[\alpha^1_u(\omega)\one_{\Lambda}
\bigl(X^\alpha (\omega)\bigr) + \alpha^2_u(
\omega)\one_{\Lambda^c}\bigl(X^\alpha(\omega )\bigr) \bigr].
\end{eqnarray*}
Since $P,P^1,P^2\in\cP$, Lemma~\ref{leadmissibilityByAlpha} yields
that $\bar{P}\in\cP$.
Moreover, we have $\bar{P}=P$ on $\cF_\tau$ and
$\bar{P}^{\tau(\omega),\omega}=P_1^{\tau(\omega),\omega}$ for
$\omega\in\Lambda$ and $\bar{P}^{\tau(\omega),\omega}=P_2^{\tau
(\omega),\omega}$ for $\omega\in\Lambda^c$. Thus $\bar{P}$ has
the required properties.
\end{pf}

We now show that the extension $\cE_t$ from Lemma
\ref{leLipschitzAndExtension} again satisfies the dynamic programming principle.
%
%th5.4 #&#
\begin{theorem}\label{thDPPesssupL1}
Let $0\leq s\leq t\leq T$ and $X\in\L^1_\cP$. Then
%
%e5.2 #&#
\begin{equation}
\label{eqDPPesssupL1} \quad\cE_s(X) = {\esssup_{P'\in\cP(s,P)}}^{
(P,\cF_s)}
E^{P'}\bigl[\cE_t(X)|\cF_s\bigr] ,\qquad P\mbox{-a.s.}\qquad
\mbox{for all } P\in\cP
\end{equation}
and, in particular,
%
%e5.3 #&#
\begin{equation}
\label{eqDPPesssupSpecialL1} \cE_s(X) = {\esssup_{P'\in\cP(s,P)}}^{
(P,\cF_s)}
E^{P'}[X|\cF_s] ,\qquad P\mbox{-a.s.}\qquad \mbox{for all } P\in\cP.
\end{equation}
\end{theorem}
\begin{pf}
Fix $P\in\cP$. Given $\eps>0$, there exists $\psi\in\UC_b(\Omega
)$ such that
\[
\bigl\|\cE_s(X)-\cE_s(\psi)\bigr\|_{L^1_\cP}\leq\|X-\psi
\|_{L^1_\cP}\leq \eps.
\]
For any $P'\in\cP(s,P)$, we also note the trivial identity
%
%e5.4 #&#
\begin{eqnarray}
\label{eqproofEsssupFormulaForL1}
&&
E^{P'}[X|\cF_s]-
\cE_s(X)
\nonumber\\
&&\qquad= E^{P'}[X-\psi|\cF_s] + \bigl(E^{P'}[\psi|
\cF_s] - \cE_s(\psi ) \bigr)\\
&&\qquad\quad{} + \bigl(\cE_s(
\psi)-\cE_s(X) \bigr) ,\qquad P\mbox{-a.s.}
\nonumber
\end{eqnarray}

(i) We first prove the inequality ``$\leq$'' in (\ref{eqDPPesssupSpecialL1}).
By (\ref{eqDPPesssupSpecial}) and Lemma~\ref{leincreasingSequence}
there exists a sequence $(P_n)$ in $\cP(s,P)$ such that
%
%e5.5 #&#
\begin{equation}
\label{eqproofEsssupFormulaForL1Psi} \cE_s(\psi)=
{\esssup_{P'\in\cP(s,P)}}^{
(P,\cF_s)} E^{P'}[\psi|\cF_s]=
\lim_{n\to\infty}E^{P_n}[\psi |\cF_s] ,\qquad P\mbox{-a.s.}
\end{equation}
Using (\ref{eqproofEsssupFormulaForL1}) with $P':=P_n$ and taking
$L^1(P)$-norms we find that
\begin{eqnarray*}
&&\bigl\|E^{P_n} [X|\cF_s]-\cE_s(X)
\bigr\|_{L^1(P)}
\\
&&\qquad\leq \|X-\psi\|_{L^1(P_n)} + \bigl\|E^{P_n}[\psi|\cF_s] -
\cE_s(\psi) \bigr\|_{L^1(P)} + \bigl\|\cE_s(\psi)-
\cE_s(X) \bigr\|_{L^1(P)}
\\
&&\qquad\leq \bigl\|E^{P_n}[\psi|\cF_s] - \cE_s(\psi)
\bigr\|_{L^1(P)} + 2\eps.
\end{eqnarray*}
Now, bounded convergence and (\ref{eqproofEsssupFormulaForL1Psi})
yield that
\[
\limsup_{n\to\infty} \bigl\|E^{P_n} [X|\cF_s]-
\cE_s(X) \bigr\|_{L^1(P)} \leq2\eps.
\]
Since $\eps>0$ was arbitrary, this implies that there is a sequence
$\tilde{P}_n\in\cP(s,P)$ such that $E^{\tilde{P}_n} [X|\cF_s]\to
\cE_s(X)$ $P$-a.s. In particular, we have proved the claimed inequality.

\mbox{}\hphantom{i}(ii) We now complete the proof of (\ref{eqDPPesssupSpecialL1}). By
Lemma~\ref{leincreasingSequence}
we can choose a sequence $(P'_n)$ in $\cP(s,P)$ such that
\[
{\esssup_{P'\in\cP(s,P)}}^{(P,\cF_s)} E^{P'}[X|\cF_s]=
\lim_{n\to\infty}E^{P'_n}[X|\cF_s] ,\qquad P\mbox{-a.s.}
\]
with an increasing limit. Let $A_n:=\{E^{P'_n}[X|\cF_s]\geq\cE_s(X)\}$.
As a result of Step (i), the sets $A_n$ increase to $\Omega$ $P$-a.s. Moreover,
\[
0\leq \bigl(E^{P'_n}[X|\cF_s]-\cE_s(X) \bigr)
\one_{A_n} \nearrow {\esssup_{P'\in\cP(s,P)}}^{(P,\cF_s)}
E^{P'}[X|\cF_s]-\cE_s(X) ,\qquad P\mbox{-a.s.}
\]
By (\ref{eqproofEsssupFormulaForL1}) with $P':=P'_n$ and by the first
equality in (\ref{eqproofEsssupFormulaForL1Psi}), we also have that
\[
E^{P'_n}[X|\cF_s]-\cE_s(X) \leq
E^{P'_n}[X-\psi|\cF_s] + \cE_s(\psi)-
\cE_s(X) ,\qquad P\mbox{-a.s.}
\]
Taking $L^1(P)$-norms and using monotone convergence, we deduce that
\begin{eqnarray*}
&&
\Bigl\|{\esssup_{P'\in\cP(s,P)}}^{(P,\cF
_s)}E^{P'}[X|
\cF_s]-\cE_s(X) \Bigr\|_{L^1(P)}
\\
&&\qquad= \lim_{n\to\infty} \bigl\| \bigl(E^{P'_n}[X|\cF_s]-
\cE_s(X) \bigr)\one_{A_n} \bigr\|_{L^1(P)}
\\
&&\qquad\leq{\limsup_{n\to\infty}} \|X-\psi\|_{L^1(P'_n)} + \bigl\|\cE_s(
\psi)-\cE_s(X)\bigr\|_{L^1(P)}
\\
&&\qquad\leq2\eps.
\end{eqnarray*}
Since $\eps>0$ was arbitrary, this proves (\ref{eqDPPesssupSpecialL1}).

(iii) It remains to show (\ref{eqDPPesssupL1}) for a given $P\in\cP
$. In view of (\ref{eqDPPesssupSpecialL1}), it suffices to prove that
\begin{eqnarray*}
&&
{\esssup_{P'\in\cP(s,P)}}^{(P,\cF_s)}  E^{P'}[X|\cF_s]
\\
&&\qquad
= {\esssup_{P'\in\cP(s,P)}}^{(P,\cF_s)} E^{P'} \Bigl[
{\esssup_{P''\in\cP(t,P')}}^{
(P',\cF_t)} E^{P''}[X|\cF_t] \big|
\cF_s \Bigr] ,\qquad P\mbox{-a.s.}
\end{eqnarray*}
The inequality ``$\leq$'' is obtained by considering $P'':=P'\in\cP
(t,P')$ on the right-hand side. To see the converse inequality, fix
$P'\in\cP(s,P)$ and choose by Lem\-ma~\ref{leincreasingSequence}
a sequence $(P''_n)$ in $\cP(t,P')$ such that
\[
{\esssup_{P''\in\cP(t,P')}}^{(P',\cF_t)} E^{P''}[X|\cF_t]=
\lim_{n\to\infty}E^{P''_n}[X|\cF_t] ,\qquad P'
\mbox{-a.s.}
\]
with an increasing limit. Then monotone convergence and the observation
that $\cP(t,P')\subseteq\cP(s,P)$ yield
\begin{eqnarray*}
&&
E^{P'} \Bigl[{\esssup_{P''\in\cP(t,P')}}^{
(P',\cF_t)}
E^{P''}[X|\cF_t] \big|\cF_s \Bigr]\\
&&\qquad =
\lim_{n\to\infty} E^{P''_n}[X|\cF_s]
\\
&&\qquad\leq {\esssup_{P'''\in\cP(s,P)}}^{(P,\cF
_s)} E^{P'''}[X|
\cF_s] ,\qquad P\mbox{-a.s.}
\end{eqnarray*}
As $P'\in\cP(s,P)$ was arbitrary, this proves the claim.
\end{pf}

We note that (\ref{eqDPPesssupSpecialL1}) determines $\cE_s(X)$ $\cP
$-q.s. and can therefore be used as an alternative definition. For
most purposes, it is not necessary to go back to the construction.
Relation (\ref{eqDPPesssupL1}) expresses the time-consistency
property of $\cE_t$. With a mild abuse of notation, it can also be
stated as
\[
\cE_s\bigl(\cE_t(X)\bigr)=\cE_s(X),\qquad 0\leq s
\leq t\leq T, X\in\L^1_\cP;
\]
indeed, the domain of $\cE_s$ has to be slightly enlarged for this
statement as in general we do not know whether $\cE_t(X)\in\L^1_\cP$.

We close by summarizing some of the basic properties of $\cE_t$.
%
%pr5.5 #&#
\begin{proposition}\label{prpropertiesOfExpectationDeterministic}
Let $X,X'\in\L^p_\cP$ for some $p\in[1,\infty)$ and let $t\in
[0,T]$. Then the following relations hold $\cP$-q.s.:
\begin{longlist}
\item $\cE_t(X)\geq\cE_t(X')$ if $X\geq X'$,
\item $\cE_t(X+X')=\cE_t(X)+X'$ if $X'$ is $\cF_t$-measurable,
\item $\cE_t(\eta X)=\eta^+\cE_t(X)+\eta^-\cE_t(-X)$ if $\eta$ is $\cF_t$-measurable and $\eta X \in\L^1_\cP$,
\item $\cE_t(X)-\cE_t(X')\leq\cE_t(X-X')$,
\item $\cE_t(X+X')=\cE_t(X)+\cE_t(X')$ if $\cE_t(-X')=-\cE_t(X')$,
\item $\|\cE_t(X)-\cE_t(X')\|_{L^p_\cP}\leq\|X-X'\|_{L^p_\cP}$.
\end{longlist}
\end{proposition}
\begin{pf}
Statements (i)--(iv) follow directly from (\ref{eqDPPesssupL1}). The
argument for (v)
is as in~\cite{Peng10}, Proposition III.2.8:
we have $\cE_t(X+X')-\cE_t(X')\leq\cE_t(X)$ by (iv) while
$\cE(X+X')\geq\cE_t(X)-\cE_t(-X')=\cE_t(X)+\cE_t(X')$ by (iv) and
the assumption on~$X'$. Of course, (vi) is contained in Lemma
\ref{leLipschitzAndExtension}.
\end{pf}

\section*{Acknowledgments}

The author thanks Shige Peng, Mete Soner and Nizar Touzi for
stimulating discussions as well as Laurent Denis, Sebastian Herrmann
and the anonymous referees for helpful comments.

%suskaldyti doi

% imsref loaded by lrinkeviciute, 2012-08-13 13:01:04
% imsref loaded by lrinkeviciute, 2012-08-13 13:11:02

\printaddresses

\end{document}